\input amstex
\input epsf
\documentstyle{amsppt}
\magnification=1200
\pretolerance=200
\tolerance=400
\TagsOnRight
\NoRunningHeads
\vcorrection{-1.2truecm}
\nologo

\define\({\left(}
\define\){\right)}

\define\C{\Bbb C}
\define\Card{\operatorname{Card}}
\define\Conv{\operatorname{Conv}}
\define\dd{\partial}
\define\diag{\operatorname{diag}}
\define\e{\operatorname{e}}
\define\Int{\operatorname{Int}}
\define\lc{\left\lceil}
\define\lpq{L_{p,q}}
\redefine\phi{\varphi}
\define\PSL{\operatorname{PSL}}
\define\R{\Bbb R}
\define\rc{\right\rceil}
\define\SL{\operatorname{SL}}
\define\sm{\setminus}
\define\VD{\operatorname{VD}}
\define\wt{\widetilde}
\define\Z{\Bbb Z}

\topmatter
\title Geometrical spines of lens manifolds	\endtitle
\author S.~Anisov \endauthor
\address Dept\. of Mathematics, Utrecht University, P.O\. Box 80{.}010, 
3508~TA Utrecht, the NETHERLANDS	\endaddress
\email anisov\@math.uu.nl \endemail
\abstract We introduce the concept of ``geometrical spine'' for 3-manifolds 
with natural metrics, in particular, for lens manifolds. We show that any 
spine of $\lpq$ that is close enough to its geometrical spine contains at least 
$E(p,q)-3$ vertices, which is exactly the conjectured value for the complexity 
$c(\lpq)$. As a byproduct, we find the minimal rotation distance (in the 
Sleator--Tarjan--Thurston sense) between a triangulation of a regular $p$-gon 
and its image under rotation.
\endabstract
\endtopmatter

\document
\head\S1. Introduction	\endhead

One can try to measure the complexity of a 3-manifold $M$ by the minimal number 
of tetrahedra $n=n(M)$ such that $M$ can be cut into $n$ tetrahedra. This 
approach has some drawbacks, in particular, this number $n$ is not necessarily 
additive under connected summation. S.~Matveev has defined the complexity 
$c(M)$ of a 3-manifold as the minimal number of vertices of an almost simple 
spine of~$M$, see~\cite6 and definitions below. The complexity $c(M)$ is equal 
to $n(M)$ for many 3-manifolds (including lenses) but is free of some of 
disadvantages of~$n(M)$. Still, it is very difficult to find $c(M)$ or even to 
estimate it from below, while reasonable upper bounds are frequently quite 
straigthforward. 

Given a metric on~$M$, one can consider the cut locus (see~\cite4) of a 
point~$x\in M$. These cut loci are spines of~$M$, though not necessarily almost 
simple, especially if the metric admits many isometries. Almost simple spines 
can, however, be obtained by small perturbations of the cut loci (or of the 
metric). This approach to constructing spines seems not to be systematically 
elaborated earlier, but it frequently yields spines with the smallest known 
numbers of vertices. We show, in particular, that this is the case for lens 
manifolds. 

Since $\lpq$ is a quotient space of $S^3\in\C^2$ by a unitary action of a 
cyclic $p$-element group~$\Z_p$ (the generator of~$\Z_p$ takes a unit vector 
$(z,w)\in\C^2$ to $(\xi z,\xi^q w)$, where $\xi=\e^{2\pi i/p}$), it carries a 
natural metric of constant curvature $+1$ and volume $2\pi^2/p$. In this paper 
we study spines of $\lpq$ that are close to cut loci of points $x\in\lpq$ with 
respect to this natural metric. 

Recall that the best known spines of $\lpq$ contain $E(p,q)-3$ vertices 
(see~\cite6; a different construction of those spines can be found in~\cite1), 
where $E(p,q)$ is the number of subtractions that the Euclid algorithm needs to 
convert an unordered pair of positive integers $(p,q)$ into $(d,0)$, where 
$d=\operatorname{g.c.d.}(p,q)$ (of course, $d=1$ if $(p,q)$ encodes a lens); in 
other words, $E(p,q)$ is equal to the sum of the elements of the continued 
fraction representing~$p/q$. 

\proclaim{Conjecture~1~\cite{6,~7}} Let $p\ge3$. Then
$$c(\lpq)=E(p,q)-3.\tag1$$
\endproclaim

In particular, for $q=1$ we have $E(p,1)=p$, and equation~(1) yields $c(L_{p,1}
)=p-3$. A possible explanation for the mysterious summand $-3$ in~(1) is that 
the number of diagonals of a $p$-gon required to triangulate it equals $p-3$, 
see also the proof of Theorem~4.

The paper is organized as follows. In \S2 we prove a number-theoretic property 
of the function $E(p,q)$ (Theorem~1), which is later used in Theorem~4. The 
proof of Theorem~1 is based on the properties of the Farey tesselation of the 
hyperbolic plane; these properties are discussed in~\S2 and in Appendix. In~\S3 
we recall the definitions of simple and special polyhedra, spines and the 
complexity of 3-manifolds. Geometrical spines are introduced in~\S4. In~\S5 we 
discuss the structure of geometrical spines of lens manifolds; in particular, 
we show the duality between spines and convex hulls (Lemma~3) and prove that 
the combinatorial type of a simple geometrical spine of a lens does not depend 
on the choice of ``base point'' involved in the construction (Lemma~4). In~\S6 
we prove (see Theorem~3) that any simple spine $P$ of a lens manifold $\lpq$ 
has at least $E(p,q)-3$ vertices provided that~$P$ is a small perturbation of a 
geometric spine. Theorem~3 is reduced to a combinatorial statement (Theorem~4) 
about the rotation distance~\cite9 between triangulations of a regular $p$-gon. 
Finally, in~\S7 we discuss sharpness of lower bounds given by Theorems~3 and~4. 

\head\S2. Digression: a sanity test for Conjecture 1	\endhead

Before going any further, let us check that $E(p,q)$ (where $(p,q)=1$ and 
$0<q<p$) is a well defined function of a lens manifold~$\lpq$. It is well known 
(see, e.g.,~\cite3) that lenses $\lpq$ and $L_{p,p-q}$ are homeomorphic, as 
well as lenses $\lpq$ and $L_{p,r}$, where $r=q^{-1}$ modulo~$p$. This suggests 
the following property of the Euclid algorithm. 

\proclaim{Theorem 1} The following relations hold\rom:
\roster
\item"a)" $E(p,q)=E(p,p-q)$\rom;
\item"b)" $E(p,q)=E(p,r)$ whenever $0<r<p$ and $rq=1$ modulo~$p$ \rom($p$ and 
$q$ are coprime in this case\rom).
\endroster
\endproclaim

If these relations were false, Conjecture~1 would be automatically false, too.

The first subtraction of the Euclid algorithm converts the pair $(p,p-q)$ into 
$(q,p-q)$. By definition of $E(p,q)$, we get $E(p,p-q)=E(q,p-q)+1$, and 
similarly $E(p,q)=E(q,p-q)+1$, which proves the first relation. The second 
relation is proved later in this section; the proof is based on some properties 
of the {\it Farey tesselation\/} of the hyperbolic plane~$H^2$.

Consider the ideal triangle in~$H^2$ with vertices at $0$, $1$, and~$\infty$.
Take its mirror images in its sides. This gives the triangles $(-1,0,\infty)$,
$(0,1/2,1)$, and~$(1,2,\infty)$, where $(a,b,c)$ denotes the ideal triangle
with vertices $a$, $b$, and~$c$. On the next step, construct the images of the
triangles obtained in the previous step under reflections in their sides that
are not sides of triangles obtained earlier. Continuing this way, we get a
tesselation of~$H^2$ into equal ideal triangles (of course, all ideal triangles 
in~$H^2$ are congruent; what is special about this tesselation is its symmetry 
in any of its  edges), see Fig.~1. It is called the {\it Farey tesselation}. 
Centers of the triangles on Fig.~1 are the vertices of the graph~$\Gamma$ dual 
to the Farey tesselation. This graph, which is the infinite binary tree embedded
in~$H^2$, is shown on~Fig.~1 by dotted lines. 

\midinsert

\epsfxsize=330pt

\centerline{\epsffile{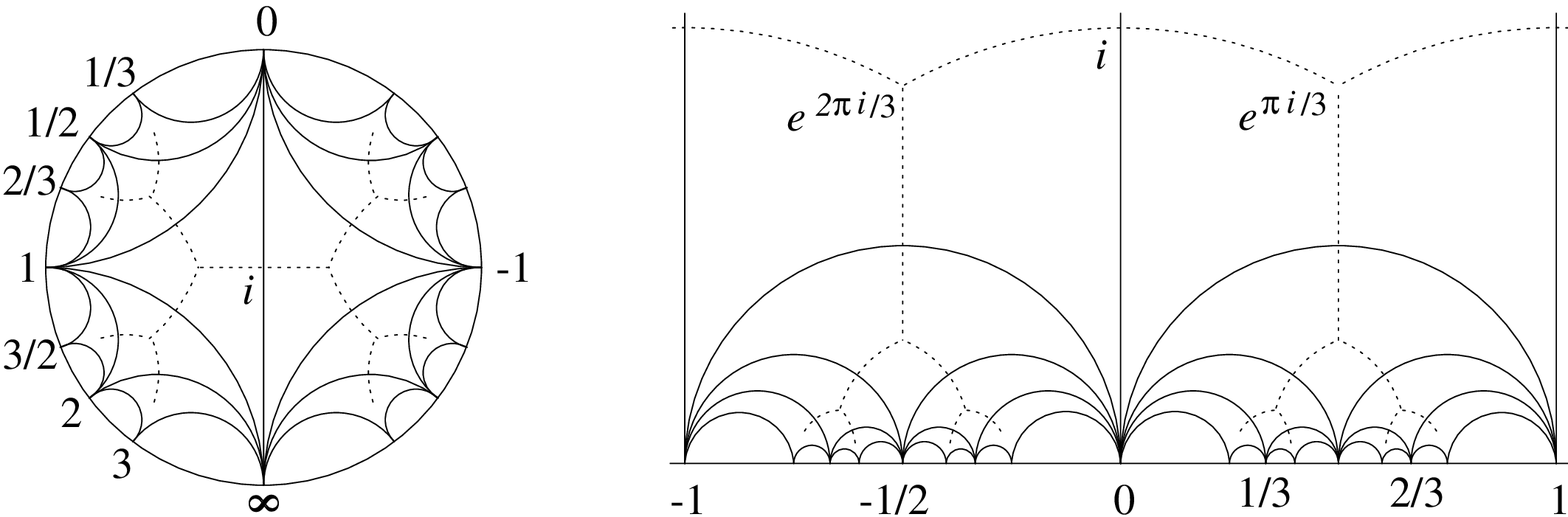}}

\botcaption{Figure 1} Farey tesselation of hyperbolic plane	\endcaption

\endinsert

Algebraic properties of the Farey tesselation are collected in Lemma~1 below; 
its geometric properties are discussed in~\S4.

\proclaim{Lemma 1} 
\roster
\item"a)" A segment $(m/n,p/q)$, where $(m,n)$ and $(p,q)$ are pairs of coprime
integers, occurs as an edge of the Farey tesselation if and only if $mq-np=
\pm1$ \rom(we agree that $\infty=1/0$\rom{);}
\item"b)" $E(p,q)$ equals the number of Farey triangles that are cut 
by the geodesic segment in~$H^2$ connecting $p/q$ with $i$ \rom(or with any 
other interior point $ti$, $t>0$, on the edge $(0,\infty)$ of the Farey 
tesselation\rom). If $|p|>|q|$, then $E(p,q)$ equals the number of Farey 
triangles cut by the geodesic segment in~$H^2$ connecting $p/q$ with $0$\rom;
\item"c)" the modular group $\PSL(2,\Z)=\SL(2,\Z)/\{\pm I\}$ acting on the upper
half-plane by fractional linear transformations $z\mapsto\frac{az+b}{cz+d}$, 
preserves the Farey tesselation\rom;
\item"d)" the third vertices of two Farey triangles incident to an edge $(m/n,
p/q)$ are $\frac{m+p}{n+q}$ and $\frac{m-p}{n-q}$\rom;
\item"e)" vertices of the Farey triangles obtained after $n$ steps of the 
construction of the tesselation and belonging to $[0,1]$ form the Farey 
sequence \rom(see~\cite{10}\rom) of depth~$n$.
\endroster
\endproclaim

For the proof of Lemma~1 see Appendix.

\demo{Proof of Theorem~\rom1} The symmetry of $H^2$ in the line $(0,\infty)$, 
$z\mapsto-\bar z$, followed by a modular mapping $z\mapsto\frac z{z+1}$, is an 
isometry of~$H^2$ preserving the Farey tesselation. This isometry is nothing but
the symmetry $z\mapsto\frac{\bar z}{\bar z-1}$ of $H^2$ in the line $(0,2)$. It 
takes $0$ to $0$ and $p/q$ to $p/(p-q)$, thus a geodesic line $(0,p/q)$ gets 
mapped to $(0,p/(p-q))$. Theorem~1a) follows now (once again) from Lemma~1b).

If $0<r<p$ and $rq=1$ modulo~$p$, we have $rq=pk+1$. The symmetry of~$H^2$ in 
$(0,\infty)$, $z\mapsto-\bar z$, followed by a modular map $z\mapsto\frac{rz+p}
{kz+q}$, is an isometry of~$H^2$ preserving the Farey tesselation. This isometry
$z\mapsto\frac{r\bar z-p}{k\bar z-q}$ takes $0$ to $p/q$ and $p/r$ to $0$, so it
takes the geodesic line $(0,p/r)$ to the geodesic $(0,p/q)$. The second 
statement of Theorem~1 follows now from Lemma~1b).\qed\enddemo

The mapping used in the proof of part b) of the theorem has a transparent 
topological meaning. Like a solid torus is obtained from $T^2\times[0,1]$ by 
contracting the meridians of $T^2\times\{0\}$, lens manifolds can be constructed
from $T^2\times[0,1]$ by contracting a family of curves of some rational 
slope~$\alpha$ in $T^2\times\{0\}$ and a family of curves of rational 
slope~$\beta$ in $T^2\times\{1\}$. Choose a basis in the lattice $\pi_1(T^2)=
\Z^2$ so that $\alpha=0$. If $\beta=p/q$ in this basis, we get~$\lpq$. To 
interchange the roles of $T^2\times\{0\}$ and $T^2\times\{1\}$ in this 
construction, we use a new basis in~$\pi_1(T^2)$ related to the previous one by 
an orientation reversing linear transformation with the matrix~$\pmatrix q&-k\\ 
p&-r\endpmatrix$. This linear map takes a vector with slope $y/x$ to a vector 
with slope $\frac{r(y/x)-p}{k(y/x)-q}$. In particular, it takes the slope $0$ to
$p/q$ and $p/r$ to $0$, thus showing that the interchanging of $T^2\times\{0\}$ 
and $T^2\times\{1\}$ in the construction above leads to encoding $\lpq$ as $L_{
p,r}$.

Alternative proofs of Theorem~1 can be found in the preprints 
\cite{1,~pp.~5--6} and~\cite{5,~pp.~29--30}. Yet another proof of Theorem~1\,b 
is based on the identity $n_k+1/(n_{k-1}+1/(n_{k-2}+\ldots+1/n_1)\dots)=p/s$ 
with $s\equiv(-1)^{k-1}r$ modulo~$p$, where $p/q=n_1+1/(n_2+1/(n_3+\ldots+1/
n_k)\dots)$ is the continued fraction representing~$p/q$. This identity holds 
whenever $n_1\ge2$ and $n_k\ge2$. 
%
%

\head\S3. Spines and complexity of 3-manifolds	\endhead

Let us recall some definitions (we follow~\cite{6,~7} here). By $K$ denote the 
1-di\-men\-si\-onal skeleton of the tetrahedron, which is just the clique (that 
is, the complete graph) with 4 vertices; $K$ is homeomorphic to a circle with 
three radii.

\definition{Definition~1} A compact 2-dimensional polyhedron is called {\it
almost simple\/} if the link of its every point can be embedded in~$K$. An
almost simple polyhedron $P$ is said to be {\it simple\/} if the link of each
point of~$P$ is homeomorphic to either a circle or a circle with a diameter or
the whole graph~$K$. A point of an almost simple polyhedron is {\it
non-singular\/} if its link is homeomorphic to a circle, it is said to be a
{\it triple point\/} if its link is homeomorphic to  a circle with a diameter,
and it is called a {\it vertex\/} if its link is homeomorphic to~$K$. The set
of singular points of a simple polyhedron~$P$ (i.e., the union of the vertices
and the triple lines) is called its {\it singular graph\/} and is denoted
by~$SP$.							\enddefinition 

It is easy to see that any compact subpolyhedron of an almost simple polyhedron
is almost simple as well. Neighborhoods of non-singular and triple points of a
simple polyhedron are shown in Fig.~2\,a,\,b; Fig.~2\,c--f represents four
equivalent ways of looking at vertices; in particular, Fig.~2\,e shows the cone 
over the 1-di\-men\-si\-onal skeleton of the tetrahedron and Fig.~2\,f shows the
Voronoi diagram of the vertices of a regular tetrahedron.

\midinsert

\epsfxsize=300pt

\centerline{\epsffile{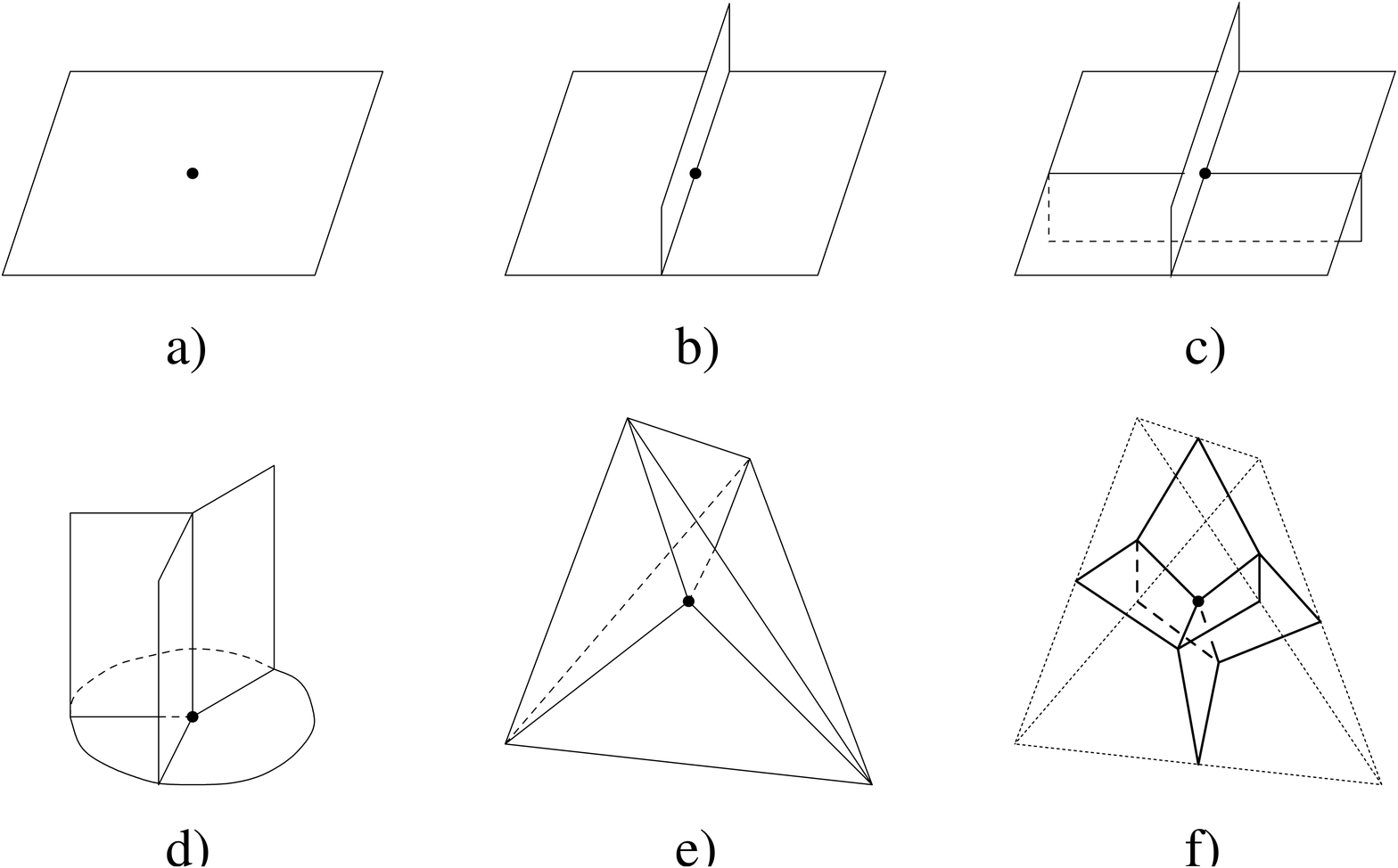}}

\botcaption{Figure 2} Nonsingular (a) and triple (b) points; ways of looking
at vertices (c--f)						\endcaption

\endinsert

\definition{Definition~2} A simple polyhedron~$P$ with at least one vertex is
said to be {\it special\/} if it contains no closed triple lines (without
vertices) and every connected component of $P\sm SP$ is a 2-dimensional cell.
								\enddefinition

\definition{Definition~3} A polyhedron $P\subset\Int M$ is called a {\it
spine\/} of a compact 3-di\-men\-sional manifold~$M$ if $M\sm P$ is
homeomorphic to $\dd M\times(0,1]$ (if $\dd M\ne0$) or to an open 3-cell (if 
$\dd M=0$). In other words, $P$ is a spine of~$M$ if a manifold $M$ with 
boundary (or a closed manifold~$M$ punctured at one point) can be collapsed 
onto~$P$. A spine $P$ of a 3-manifold $M$ is said to be {\it almost simple}, 
{\it simple}, or {\it special\/} if it is an almost simple, simple, or special 
polyhedron, respectively.					\enddefinition

Given a special spine $P$ of a compact manifold $M^3$, one can construct a dual 
singular triangulation of $M^3$ with one vertex (lying in the middle of the 
3-cell $M\sm P$), see Fig.~2\,f; if $M$ is a manifold with connected boundary, 
a similar construction gives a triangulation of the one-point compactification 
of $M\sm\dd M$. In both cases, there is a one-to-one correspondence between 
vertices of~$P$ and tetrahedra of the triangulation.

\definition{Definition~4} The {\it complexity\/} $c(M)$ of a compact
3-manifold~$M$ is the minimal possible number of vertices of an almost simple
spine of~$M$. An almost simple spine with the smallest possible number of
vertices is said to be a {\it minimal\/} spine of~$M$.		\enddefinition

\proclaim{Theorem~2~\cite6} Let $M$ be an orientable irreducible $3$-manifold 
with incompressible \rom(or empty\rom) boundary and without essential annuli. 
If $c(M)>0$ \rom(that is, if $M$ is different from \rom(possibly punctured\rom) 
$S^3$, $\R P^3$, and $L_{3,1}$\rom), then any minimal spine of~$M$ is special. 
								\endproclaim

Note that lenses $\lpq$ with $p>3$ satisfy the conditions of Theorem~2. Thus any
minimal spine of $\lpq$ is special and corresponds to a decomposition of~$\lpq$
into $n=c(\lpq)$ tetrahedra.

\head\S4. Cut loci and Voronoi diagrams \endhead

We are going to use relations between spines, cut loci, and Voronoi diagrams.
For the definition and some properties of cut loci, see \cite{4, Chapter~VIII, 
\S7}. Here we only recall that given a point $x\in M$ of a complete Riemannian 
manifold  $(M,g)$, the {\it cut locus\/} of $x$ is the closure (in~$M$) of the 
set of points $y\in M$ such that the shortest geodesic between $x$ and $y$ is 
not unique. If $M$ is compact, then the preimage $\exp^{-1}(C(x))\subset T_xM$ 
of a cut locus $C(x)$ under the exponential map $\exp\:T_xM\to M$ is 
homeomorphic to a sphere and $M\sm C(x)$ is homeomorphic to a ball~\cite4.

Now we have the following method of constructing spines of compact 
$3$-mani\-folds: choose a metric $g$ on $M$ and a point~$x\in M$. Then the cut 
locus $C(x)$ is a spine of~$M$, though not necessarily special. Contrary, for 
any spine $P$ of $M$, there is a metric~$g$ such that $P$ is isotopic to the 
cut locus $C(x)$ for some $x\in M$ with respect to~$g$ (take a standard metric 
on a unit ball homeomorphic to $M\sm P$, transfer it to $M$ and smoothen it 
around~$P$).

\example{Examples} 1.~The cut locus of the sphere~$S^n$ (with the standard 
metric) with respect to the North pole~$N$ is the South pole~$S$. So a single 
point is a spine of~$S^3$.

2.~The cut locus of a ``rectangular'' flat torus $T^2$ with respect to a point 
$x=(\phi_0,\psi_0)$ is the union of the parallel $\psi=\psi_0+\pi$ and the 
meridian $\phi=\phi_0+\pi$ opposite to~$x$. The complement $T^2\sm C(x)$ is a 
flat rectangle. 

3.~The cut locus of a generic flat torus~$T^2$ with respect to any point $x\in 
T^2$ is the union of three nonhomotopic geodesic segments connecting the two 
local maxima of~$d_x$, the distance to~$x$. The complement $T^2\sm C(x)$ is a 
flat centrally symmetric hexagon. The preimage of $C(x)$ under the universal 
covering $\R^2\to T^2$ is shown on Fig.~3. 

4.~The cut locus of the torus $T^3=\R^3/\Z^3$ (where $\R^3$ carries a standard 
metric and $\Z^3$ is the integer span of an orthogonal basis) with respect to 
the point $(1/2,1/2,1/2)$ is the union of three ``coordinate'' tori $T^2$. This 
spine is not special, since it contains, in particular, lines of transversal 
intersection of two surfaces. However, the cut locus of a point in a generic 
flat $T^3$ is a special spine with $6$ vertices, which is a minimal spine, 
because $c(T^3)=6$~\cite7. See also Fig.~12 in~\cite1.		\endexample

{\it Voronoi diagrams\/} are, roughly speaking, cut loci with respect to many 
points. Originally they were defined in~\cite{11} for Euclidean plane (or, more 
generally, for $\R^n$) with finitely many nodes $A_1,\dots,A_k\in\R^n$ as the 
set of points where the nearest node is not unique. A Voronoi diagram 
divides~$\R^n$ into Voronoi domains $U_1,\dots,U_k$, where 
$U_i$ consists of the points of $\R^n$ that are closer to $A_i$ than to $A_j$ 
with any $j\ne i$. This admits an obvious generalization for locally finite 
subsets of complete Riemannian manifolds. The following examples 5--8 are 
parallel to the examples 1--4 above.

\example{Examples} 5.~The Voronoi diagram of the sphere $S^n$ with respect to 
a single point~$N$ is empty; the whole $S^n$ forms the single Voronoi domain. 
The Voronoi diagram of $S^n$ with respect to two antipodal points $N$ and $S$ 
(North and South poles) is the equator $S^{n-1}\subset S^n$; Northern and 
Southern hemispheres form the Voronoi domains (of course, we use here and below 
the standard metric on~$S^n$).

6.~The Voronoi diagram of the square lattice $\{(m,n)\mid m,n\in\Z\}$ in~$\R^2$ 
is the square grid formed by lines $\{x\}=1/2$ and $\{y\}=1/2$, where $\{\cdot
\}$ denotes the fractional part. The Voronoi domains are unit squares centered 
at points with integer coordinates.

7.~The Voronoi diagram of the lattice $\{(m+n/2,\sqrt3n/2)\mid m,n\in\Z\}
\subset\R^2$ generated by the vectors $(1,0)$ and $(1/2,\sqrt3/2)$ looks like a 
honeycomb, see Fig.~3, left. The Voronoi domains are regular hexagons. The 
Voronoi diagram of a generic lattice in~$\R^2$ is also hexagonal, see Fig.~3,
right.

\midinsert

\epsfxsize=300pt

\centerline{\epsffile{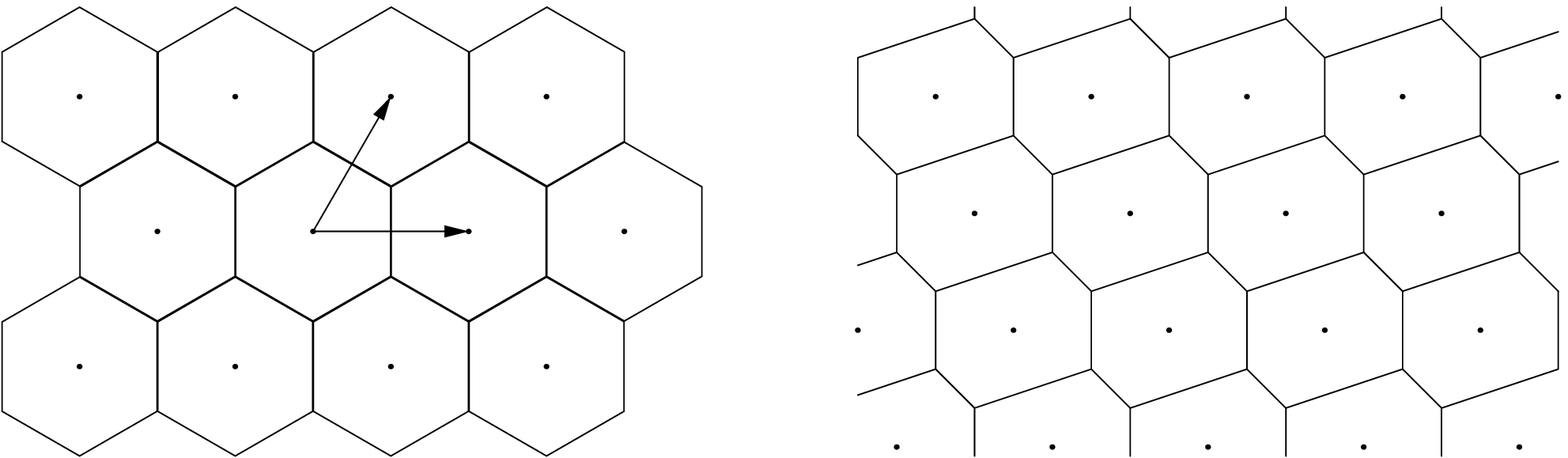}}

\botcaption{Figure 3} Voronoi diagrams of generic lattices in~$\R^2$ 
\endcaption

\endinsert
\endexample

8.~The Voronoi diagram of the cubic lattice $\{(k,l,m)\mid k,l,m\in\Z\}$ in~$\R
^3$ is formed by the planes $\{x\}=1/2$, $\{y\}=1/2$, and $\{z\}=1/2$. The 
Voronoi domains are unit cubes. For a generic rank~3 lattice in~$\R^3$, the 
Voronoi domains are centrally symmetric polyhedra with 14 facets, 36 edges, and 
24 vertices.

The next two examples are central for the rest of the paper.

\example{Example 9: universal covering} Let $p\:\wt M\to M$ be the universal 
covering of~$M$, where $\wt M$ is a constant curvature space $H^n$, $\R^n$, 
or~$S^n$ (but $M$ itself is different from $S^n$); assume that $M$ inherits a 
constant curvature metric from~$\wt M$. Choose $x\in M$, set $X=\{x_1,x_2,\dots
\}=p^{-1}(x)\subset\wt M$, and consider the Voronoi diagram in $\wt M$ with 
respect to~$X$. In this situation all Voronoi domains are contractible and the 
projection of the Voronoi diagram in~$\wt M$ is the cut locus of~$M$ with 
respect to~$x$. If $M$ is $\R P^3$, we get Example~5. Other examples above also 
are particular cases of this construction.			\endexample

\example{Example 10: Farey tesselation} The Teichm\"uller space for~$T^2$ is 
the hyperbolic plane~$H^2=\{z=x+iy\in\C\mid y>0\}$: $T^2_z$ can be thought of 
as the quotient space of $\R^2$ over the lattice $\{m\cdot1+n\cdot z\mid m,n\in
\Z\}\subset\C$. Let $X\subset H^2$ be the set of all parameters $z$ 
corresponding to the tori with three equally short shortest geodesics (i.e., 
tori glued from a regular hexagon). Then the Farey tesselation is nothing but 
the Voronoi diagram of $H^2$ with respect to~$X$. The vertices of the triangles 
represent the slopes of the elements of $\pi_1(T^2)$, and the modular group 
action on $H^2$ (restricted to the absolute) corresponds to coordinate changes 
in $\pi_1(T^2)$. The intersection points of solid and dashed lines on Fig.~1 
represent ``square'' tori, and other points on the solid lines represent 
``rectangular'' tori. The cut locus $\theta$ of a flat torus (generically it is 
a graph with two vertices and three pairwise nonhomotopic edges) changes 
isotopically as its parameter $z$ varies inside of a Farey triangle, and 
bifurcates when $z$ crosses Farey edges. Three vertices of a Farey triangle 
containing $z$ are the slopes of three cycles in $T^2_z$ formed by pairs of 
edges of~$\theta$. We omit rather straightforward proofs of these statements; 
see also~\cite2.						\endexample

\head\S5. Application to lens manifolds \endhead

According to Example~9 above, a spine of a lens manifold $\lpq$ can be obtained 
as the covering projection image of the Voronoi diagram (in $S^3\subset\R^4$) 
of the $\Z_p$-orbit $\Z_p(z,w)=\{(\xi^kz,\xi^{kq}w)\mid k=0,1,\dots,p-1,\,\xi=
\e^{2\pi i/p}\}$ of a point $(z,w)\in\C^2$. Then the ambient space $\R^4$ is 
divided into $p$ congruent convex polyhedral cones with a common vertex at the 
origin. The Voronoi domains in $S^3$ are described by the following statement.

\proclaim{Lemma 2} Let $\{A_1,\dots,A_k\}\subset S^{n-1}$ be a finite set of 
points on the unit sphere in~$\R^n$. This set defines Voronoi domains in $S^{n-
1}$ and in $\R^n$. Then the Voronoi domains in $S^{n-1}$ are the intersections 
of the Voronoi domains in $\R^n$ with~$S^n$. 			\endproclaim

\demo{Proof} This is obvious for a two-point set $\{A_1,A_2\}$. The general 
case follows from the case of two points. \qed\enddemo

\proclaim{Lemma 3} Let $A=\{A_1,\dots,A_k\}\subset S^{n-1}$ be a finite set of 
points on the unit sphere in~$\R^n$. Suppose that this set is not contained in 
a hyperplane. By $\Conv(A)$ denote the convex hull of $A\subset\R^n$. Then the 
Voronoi diagram $\VD(A)\subset S^{n-1}$ is dual to~$\Conv(A)$\rom: its vertices 
are unit outer normals to the facets of \,$\Conv(A)$ etc. In particular, the 
combinatorial type of \,$\VD(A)$ is determined by the combinatorial type 
of \,$\Conv(A)$.		\endproclaim

\demo{Proof} This follows from the definitions and from Lemma~2. \qed\enddemo

Consider the circles $\wt S^1_z=\{(\e^{i\phi},0)\mid0\le\phi<2\pi\}$ and $\wt 
S^1_w=\{(0,\e^{i\phi})\mid0\le\phi<2\pi\}$ in~$S^3$. They cover two {\it core 
circles\/} $S^1_z$ and $S^1_w$ in~$\lpq$. Take a point $x=(1,0)$ (or any other 
point of $\wt S^1_z$). Its orbit $\Z_px$ is a regular $p$-gon in the plane $w=
0$. The Voronoi diagram for $\Z_px$ in $\R^4$ looks like ``an orange slice 
times $\R^2$'' and consists of $p$ copies of $\R^3_+$ glued together along the 
plane $z=0$ so that all dihedral angles are equal to~$2\pi/p$. The Voronoi 
diagram for $\Z_px$ in $S^3$ consists of $p$ hemispheres $S^2_+=S^2\cap\R^3_+$ 
glued together at equal angles along the circle~$\wt S^1_w$.

The covering mapping takes each of these hemispheres to a disk so that its 
boundary goes $p$ times along~$S^1_w$. The resulting spine of $\lpq$ consists 
of this disk and~$S^1_w$. Locally, its transversal (to~$S^1_w$) section looks 
like the set $Y_p=\{r\e^{2\pi k i/p}\mid0\le r<1,\,k=0,1,\dots,p-1\}$ ($Y_3$ 
looks like Y, $Y_4$ like $\times$, $Y_5$ like $\star$, $Y_6$ like $\ast$ etc.). 
A neighborhood of $S^1_w$ in this spine fibers over $S^1_w$ with fiber~$Y_p$.
This fiber bundle is nontrivial: monodromy is a positive (clockwise as the 
parameter $\phi$ on~$S^1_w$ grows) rotation by $2\pi r/p$, where $r=q^{-1}$ 
modulo~$p$.\footnote{The same construction with a point $y=(0,1)$ instead of 
$x=(1,0)$ gives a $Y_p$-fibration over~$S^1_z$ with rotation by $2\pi q/p$ as 
the monodromy.} 

If $p=3$, this construction yields a simple spine of $L_{3,1}$ without vertices
(thus showing that $c(L_{3,1})=0$). However, for $p>3$ the spines obtained this 
way are not almost simple because of the line (the core circle~$S^1_w$) of 
multiplicity~$p$. Simple spines can be obtained by small perturbation of these 
cut loci or, in many cases, by choosing the point $x$ outside of the core 
circles. 

\proclaim{Lemma 4} If a cut locus $C(x)$ is a simple spine of $\lpq$, then its 
combinatorial type is independent of the choice of $x$ in $\lpq\sm\{S^1_w,S^1_z
\}$.								\endproclaim

\demo{Proof} Suppose that $C(x)$ is a simple spine of $\lpq$. Let $(z_0,w_0)$ 
be one of $p$ preimages of $x$ under the covering $S^3\to\lpq$. Then the set of 
all preimages of $x$ is the orbit $\Z_p(z_0,w_0)$, and its Voronoi diagram in 
$S^3$ is a simple polyhedron, because it covers~$C(x)$. By Lemma~3, the 
combinatorial type of~$C(x)$ is determined by the combinatorial type of the 
convex hull $\Conv(\Z_p(z_0,w_0))$. Note that $z_0\ne0$ and $w_0\ne0$; 
otherwise $C(x)$ is not a simple polyhedron.

Let $(z_1,w_1)$ be any point of $S^3\sm\{\wt S^1_w,\wt S^1_z\}$, i.e., $z_1^2+w
_1^2=1$, $z_1\ne0$, $w_1\ne0$. Consider a linear transformation of $\C^2$ 
defined by a diagonal matrix $\diag(z_1/z_0,w_1/w_0)$. This is an invertible 
linear transformation that takes the orbit $\Z_p(z_0,w_0)$ to $\Z_p(z_1,w_1)$. 
It also takes $\Conv(\Z_p(z_0,w_0))$ to $\Conv(\Z_p(z_1,w_1))$. Then these two 
convex hulls have the same combinatorial type. Therefore, the combinatorial 
type of $C(x)$ does not depend on the choice of $x\in\lpq\sm\{S^1_w,S^1_z\}$.
\qed\enddemo

If $q=\pm1$ modulo~$p$, then the $\Z_p$-orbit of any point is a (flat) regular 
$p$-gon. In this case $C(x)$ cannot be a simple spine of~$\lpq$ unless $p=3$. 
If $q\ne\pm1$ modulo~$p$, then $C(x)$, where $x\in\lpq\sm\{S^1_w,S^1_z\}$, is 
a simple spine of $\lpq$ with $E(p,q)-3$ vertices, so the assumption of 
simplicity of~$C(x)\subset\lpq$ is satisfied, see~\cite2.

\head\S 6. Perturbations of $C(x)$ and spines of $\lpq$ \endhead

Recall (see \S4) that any spine $P\subset M^3$ is (up to isotopy) a cut locus 
$C(x)$ for some pair $(x,g)$, where $x\in M$ and $g$ is a Riemanian metric. We 
say that a spine $P'$ is a {\it small perturbation\/} of $P$ if $P'$ is a cut 
locus for a pair $(x',g')$ that is a small perturbation of~$(x,g)$.

\proclaim{Theorem~3} If an almost simple spine $P$ of $\lpq$ is a small 
perturbation of a cut locus $C(x)$ \rom(defined by the standard metric on 
$\lpq$\rom), then $P$ contains at least $E(p,q)-3$ vertices. 	\endproclaim

\demo{Proof} It is sufficient to prove the theorem for the case $x\in S^1_z$. 
Indeed, the case $x\in S^1_w$ is similar to the case $x\in S^1_z$, and if $x'$ 
does not lie on a core circle, the combinatorial type of $C(x')$ (and of its 
small perturbations) does not depend on $x'$, so we can assume that $x'$ lies 
very close to $S^1_z$; then $C(x')$ is a small perturbation of $C(x)$ with 
$x\in S^1_z$, and small perturbations of $C(x')$ also are small perturbations 
of~$C(x)$.

A $C^\infty$-small perturbation of the metric leads to a small perturbation 
of~$C(x)$. This perturbation is a local isotopy in a neighborhood of any 
nonsingular point of $C(x)$, of any point on a triple line of $C(x)$, and in a 
neighborhood of any vertex, see Fig.~1. (Thus a small perturbation of a simple 
polyhedron is a simple polyhedron of the same combinatorial type, in 
particular, with the same number of vertices.)

\midinsert

\epsfxsize=240pt

\centerline{\epsffile{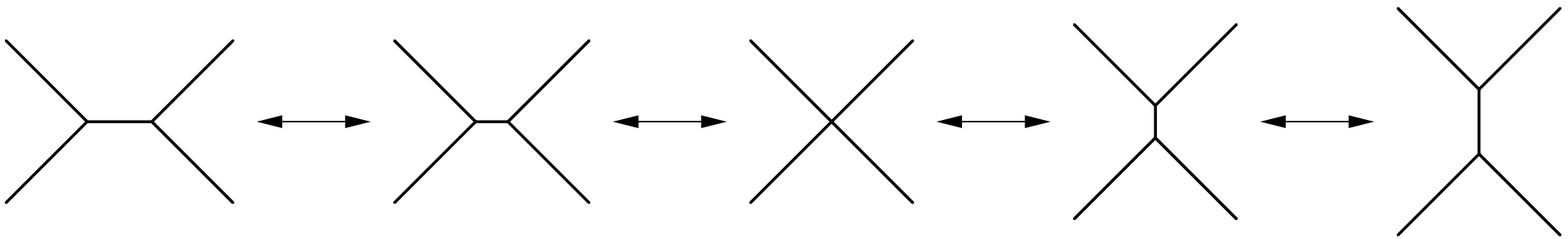}}

\botcaption{Figure 4} A flip 
\endcaption

\endinsert

Consider a cut locus $C(x)\subset\lpq$, where $x\in S^1_z$; it was described in 
the paragraph following Lemma~3. If a perturbation is generic, multiple line 
(in our case $S^1_w$) splits in a number of triple lines, which may end in new 
vertices. Consider a section of $C(x)$ by a transversal to $S^1_w$ at $y\in S^1
_w$. For a generic $y$, the effect of the perturbation is the splitting of a 
degree $p$ vertex of the graph~$Y_p$ into $p-2$ trivalent vertices; denote the 
perturbed section by~$Y_p(y)$. As $y\in S^1_w$ varies, {\it flips\/} (see 
Fig.~4) may occur in~$Y_p(y)$. Note that flips occuring in transversal sections 
of a simple polyhedron $P$ correspond to vertices of~$P$; see Fig.~5, where 
transversal sections of a neighborhood of a vertex in a simple polyhedron by 
three parallel planes are shown.

\midinsert

\epsfxsize=210pt

\centerline{\epsffile{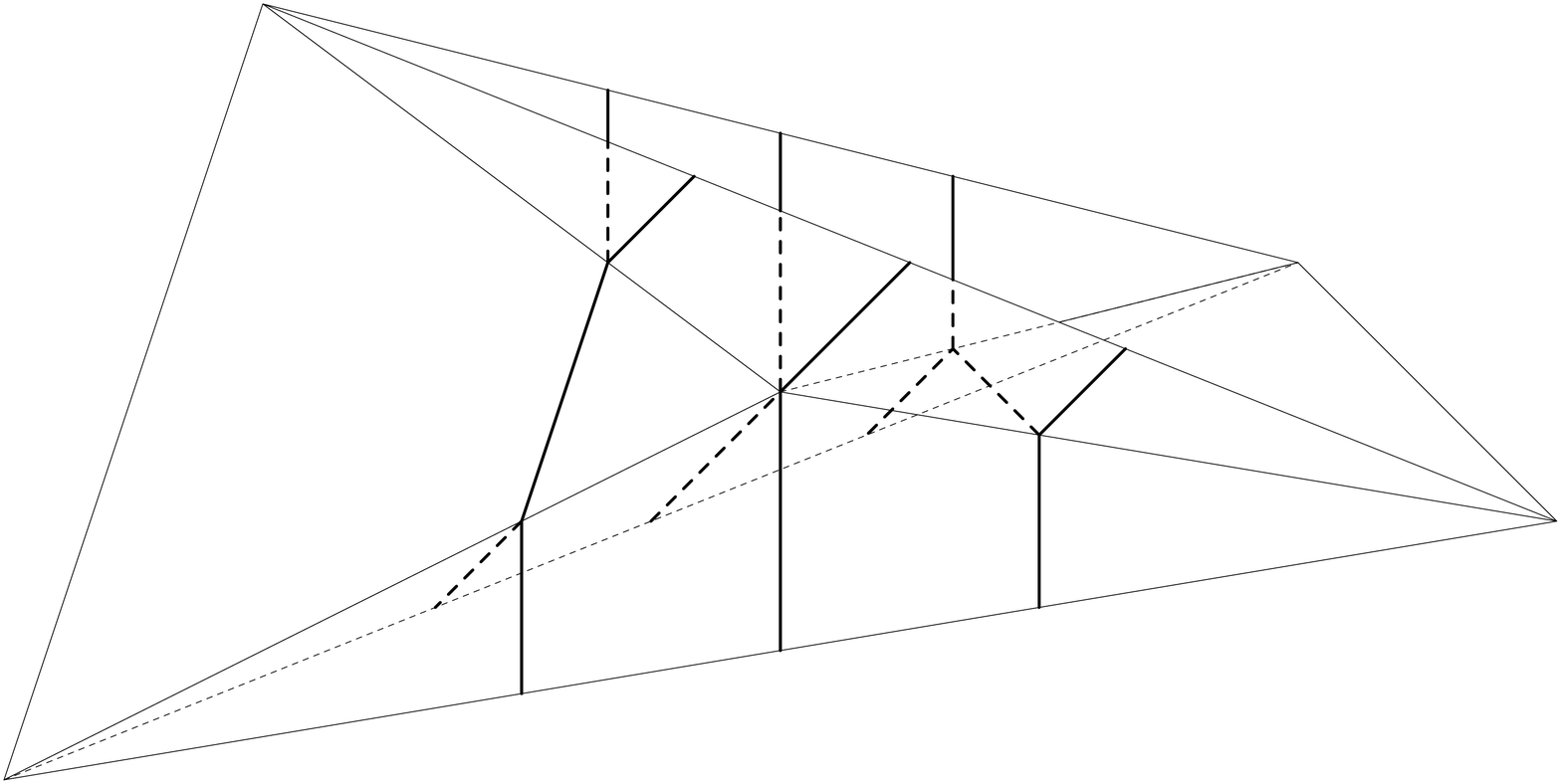}}

\botcaption{Figure 5} Vertices correspond to flips 
\endcaption

\endinsert

Let us move $y$ along the circle $S^1_w$. Then $Y_p(y)$ undergoes isotopy and 
flips, and the sequence of flips converts $Y_p(0)$ into $Y_p(2\pi)$; due to 
monodromy (described in~\S5), the latter graph is the former one rotated by the 
angle $2\pi r/p$, where $r=q^{-1}$ modulo~$p$. We have to estimate from below 
the number of flips required to convert $Y_p(0)$ into its image under rotation. 

There is a natural one-to-one correspondence between the isotopy classes of 
trivalent resolutions of $Y_p$ (with fixed $p$ boundary points) and the 
triangulations of the regular $p$-gon: for any triangulation, its dual graph is 
an isotopy class of trivalent resolutions of $Y_p$, and vise versa. A {\it 
flip\/} in a triangulation is replacing of two triangles $ABC$ and $ACD$ having 
a common side $AC$ by two triangles $BCD$ and $ABD$ with a common side $BD$; in 
other words, a flip switches the diagonal in a triangulated quadrilateral. 
Flips in trivalent resolutions of $Y_p$ correspond to flips of dual 
triangulations. Taking into account Theorem~1\,b, Theorem~3 follows from 
Theorem~4 below. \qed\enddemo

\proclaim{Theorem 4} Let $\Delta_1$, $\Delta_2$ be two triangulations of a 
regular $p$-gon, $p\ge3$. By $d(\Delta_1,\Delta_2)$ denote the rotation 
distance \rom(see~\rom{\cite9)} between $\Delta_1$ and~$\Delta_2$, that is, 
the minimal number of flips required to convert $\Delta_1$ into $\Delta_2$. If 
$\Delta_2$ is $\Delta_1$ rotated by $2\pi q/p$, where $p>q>0$ and $(p,q)=1$, 
then $d(\Delta_1,\Delta_2)\ge E(p,q)-3$.		\endproclaim

\remark{Remarks} 1. The lower bound given by Theorem~4 is exact, i.e., there 
always exists a triangulation $\Delta_1$ such that $d(\Delta_1,\Delta_2)=
E(p,q)-3$.  \newline
2. The condition $(p,q)=1$ is not necessary for the inequality $d(\Delta_1,
\Delta_2)\ge E(p,q)-3$. The equality $d(\Delta_1,\Delta_2)=\max(0,E(p,q)-3)$ 
holds for all $p$, $q$ such that $p\ge q\ge0$ and $p\ge3$. Nevertheless, in the 
proof below we assume, for the sake of simplicity, that $(p,q)=1$. \endremark

\demo{Proof} Without loss of generality, we may suppose that $p>2q$. Indeed, 
clockwise rotation by $2\pi q/p$ is the same as counterclockwise rotation by $2
\pi(p-q)/p$, and $E(p,q)=E(p,p-q)$ by Theorem~1\,a. The number $E(p,q)$ appears 
in this proof as the sum of the elements of the continued fraction expansion 
of~$p/q$. 

We start with the case $q=1$; then $E(p,q)=E(p,1)=p$. A triangulation of a 
$p$-gon involves $p-3$ diagonals. Any of them intersects its own image under 
the ``minimal'' rotation by $2\pi/p$, therefore every of them needs to be moved 
by a flip, and this requires at least as many flips as diagonals, that is, at 
least $p-3$ flips. This proves the theorem for $q=1$, i.e., for the case of the 
shortest continued fraction. 

Two simple ideas have been used here: first, $d(\Delta_1,\Delta_2)$ is not less 
than the number of diagonals of $\Delta_1$ that are not mapped to other 
diagonals of $\Delta_1$ by the rotation, so it is sufficient to estimate the 
number of those diagonals, and, second, long diagonals cannot survive small 
rotations. 

Define the {\it length\/} of a diagonal to be the number of sides of the 
polygon in the shorter arc bounded by the endpoints of the diagonal; thus, the 
length of a diagonal is at least 2 and at most $p/2$. Note that any diagonal 
longer than~$x$ intersects its image under a rotation by $2\pi x/p$. \enddemo

\proclaim{Lemma 5} Let $x$ be an integer such that $2\le x\le p/2$. Then any 
triangulation of a regular $p$-gon contains at least $\lc p/x\rc-3$ diagonals 
that are longer than~$x$. 		\endproclaim

\demo{Proof} Cut the polygon along every short diagonal of the triangulation 
(that is, along every diagonal of length at most~$x$). Consider the piece of 
the $p$-gon containing its center. This piece has at least $\lc p/x\rc$ sides 
and is triangulated by at least $\lc p/x\rc-3$ diagonals, all of which are 
longer than~$x$. 						\qed\enddemo

Now consider the case of a two-term continued fraction, $p/q=n_1+1/n_2$. Then 
$p=n_1n_2+1$, $q=n_2$, $E(p,q)=n_1+n_2$, and the number of diagonals in a 
triangulation is $p-3=n_1n_2-2$. A (clockwise) rotation by $2\pi q/p$, repeated 
$n_1$ times, is equivalent to a (counterclockwise) rotation by $2\pi/p$, which 
none of the $p-3$ diagonals can survive. Therefore, at least $\lc\frac{p-3}
{n_1}\rc$ diagonals are destroyed by a single rotation by~$2\pi q/p$; this 
number can be equal to $n_2$ or $n_2-1$. However, one can give a better 
estimate. 

By Lemma~5, the triangulation contains at least $\lc p/q\rc-3=n_1-2$ diagonals 
that are longer than~$q$. A rotation by~$2\pi q/p$ destroys all of them and at 
least $1/n_1$ fraction of the other diagonals (because otherwise some of the 
latter diagonals survive a rotation by~$2\pi q/p$ repeated $n_1$ times), and we 
get $d(\Delta_1,\Delta_2)\ge(n_1-2)+\lc\frac{(p-3)-(n_1-2)}{n_1}\rc=(n_1-2)+
\lc\frac{n_1n_2-n_1}{n_1}\rc=(n_1-2)+(n_2-1)=E(p,q)-3$.

In the general case we have $p/q=n_1+1/(n_2+1/(n_3+\ldots+1/n_k)\cdots)$. The 
$i$th {\it convergent\/} of this continued fraction is $p_i/q_i=n_1+1/(n_2+
\ldots+1/n_i)\cdots)$, $1\le i\le k$; we also set $p_0=1$ and $q_0=0$. Further, 
set $r_0=p$, $r_1=q$ and define the numbers $r_2,\dots,r_k$ as the remainders 
in the Euclid algorithm: the relations $r_{i-1}=n_ir_i+r_{i+1}$ define the 
$r_i$ recurrently. Note that $r_k=(p,q)=1$. It is convenient to set~$r_{k+1}=0$.

\proclaim{Lemma 6} For $i=1,2,\dots,k-1$, we have $p/r_i>p_i$.
\endproclaim

\demo{Proof} This follows from the relation $p=p_ir_i+p_{i-1}r_{i+1}$, where 
$i=1,\ldots,k-1$. This relation is proved by induction on~$i$. For $i=1$ it 
takes the form $p=p_1r_1+p_0r_2=n_1q+r_2$, which is the first step of the 
Euclid algorithm. The induction step is the equation $p_ir_i+p_{i-1}r_{i+1}=
p_{i+1}r_{i+1}+p_ir_{i+2}$. The relations $p_{i+1}=p_in_{i+1}+p_{i-1}$ and $r_i
=n_{i+1}r_{i+1}+r_{i+2}$ show that both sides are equal to $n_{i+1}p_ir_{i+1}+
p_{i-1}r_{i+1}+p_ir_{i+2}$. 					\qed\enddemo

Consider a rotation by $2\pi q/p$ repeated $p_i$ times, $i=0,1,\dots,k$. The 
following lemma shows that this gives a (positive or negative) rotation by 
$2\pi r_{i+1}/p$.

\proclaim{Lemma 7} For $i=0,1,\dots,k$, the congruence $p_iq\equiv(-1)^ir_{i+1}
\mod p$ holds. 							\endproclaim

\demo{Proof} Induction on~$i$. The base, $i=0$, is obvious (recall that $p_0=1$ 
and $r_1=q$). Induction step: $p_{i+1}q=(n_{i+1}p_i+p_{i-1})q\equiv(-1)^in_{i+1}
r_{i+1}+(-1)^{i-1}r_i=(-1)^{i+1}r_{i+2}\mod p$, because $r_{i+2}=r_i-n_{i+1}r_
{i+1}$. 							\qed\enddemo

\proclaim{Lemma 8} If $2\le j\le k$, then the following relation holds\rom:
$$(n_1-2)+\sum_{i=2}^{j-1}p_{i-1}n_i+p_{j-1}(n_j-1)=p_j-3.\tag2$$
\endproclaim

\demo{Proof} Since $p_0=1$, equation (2) is equivalent to $\sum_{i=1}^jp_{i-1}
n_i=p_j+p_{j-1}-1$. This is easy to prove by induction, because $p_i=n_ip_{i-1}
+p_{i-2}$, see, e.g.,~\cite{10}.				\qed\enddemo

Let us return to the proof of Theorem~4. By $S_i$, $i=1,\dots,k$, denote the 
set of diagonals of~$\Delta_1$ that are longer than $r_i$ but not longer than 
$r_{i-1}$; set $s_i=\Card(S_i)$. Lemma~5 implies that $s_1+\ldots+s_i\ge
\lc p/r_i\rc-3$; therefore, by Lemma~6,
$$\align
s_1+\ldots+&s_i>p_i-3\qquad\text{for\ \ } i=1,\dots,k-1\tag3{${}_i$}\\
s_1+\ldots+&s_k=p-3.\tag3
\endalign$$

\proclaim{Lemma~9} Let $s_1,\dots,s_k$ be integers. Let $j$ be an integer such 
that $2\le j\le k$.
\roster
\item"a)" Suppose that $(3_i)$ holds for any $i=1,\dots,j-1$, and $s_1+\ldots+
s_j=p_j-3$. Then 
$$s_1+\lc s_2/p_1\rc+\ldots+\lc s_j/p_{j-1}\rc\ge 
	n_1+\ldots+n_j-3.\tag4{${}_j$}$$
\item"b)" Under the assumptions of \rom{a)},  the inequality $(4_j)$ is an 
equality if and only if $s_1=n_1-2$, $s_l=p_{l-1}n_l$ for $l=2,\dots,j-1$, and 
$s_j=p_{j-1}(n_j-1)$.
\item"c)" Suppose that $(3_i)$ holds for any $i=1,\dots,j-1$, and $s_1+\ldots+
s_j\ge p_j-3+dp_{j-1}$, where $d$ is a positive integer. Then 
$$s_1+\lc s_2/p_1\rc+\ldots+\lc s_j/p_{j-1}\rc\ge 
	n_1+\ldots+n_j-3+d.\tag5{${}_j$}$$
\item"d)" Under the assumptions of \rom{c)},  the inequality $(5_j)$ is an 
equality if and only if $s_1=n_1-2$, $s_l=p_{l-1}n_l$ for $l=2,\dots,j-1$, and 
$s_j=p_{j-1}(n_j-1+d)$.
\endroster
\endproclaim

\demo{Proof} If $s_1=n_1-2$, $s_l=p_{l-1}n_l$ for $l=2,\dots,j-1$, and $s_j=
p_{j-1}(n_j-1)$, then the assumptions of the part a) are satisfied by virtue of 
Lemma~8. The ``if'' part of the statement~b) is obvious. Other statements are 
proved by induction on~$j$. The idea is as follows: given the sum of the~$s_i$, 
the left hand side of~(4) and~(5) would be minimized if, roughly speaking, as 
much as possible were divided by the biggest denominator (that is, by the last 
one), but the constraints~(3${}_i$) and~(3) do not allow us to overload the 
last summand~$s_j$; the ceiling function $\lc\cdot\rc$ is responsible for the 
clause ``roughly speaking'' in the above. 

Base of induction: $j=2$. Under the assumptions of~a), if $s_2=p_1(n_2-1)$, 
then $s_1=(s_1+
s_2)-s_2=p_2-3-s_2=n_2p_1+p_0-3-p_1(n_2-1)=n_1-2=p_1-2$ and $s_1+\lc s_2/p_1\rc
=n_1+n_2-3$. If $s_2>p_1(n_2-1)$, then $s_1=p_2-s_2<p_1-2$, which 
contradicts~($3_1$). If $s_2<p_1(n_2-1)$, then $s_2=p_1(n_2-1)-t$ and $s_1=n_1-
2+t$ for some $t>0$. Then $s_1+\lc s_2/p_1\rc=n_1+n_2-3+t+\lc-t/p_1\rc>n_1+n_2-
3$, because $p_1=n_1\ge2$ (due to the assumption $p>2q$). This proves 
statements a) and b) for~$j=2$. The proof of the statements c) and d) for $j=2$ 
is similar.

Induction step: from $j=l-1$ to $j=l$. Under the assumptions of~a), if $s_l=p_
{l-1}(n_l-1)$, then $s_1+\ldots+s_{l-1}=p_l-3-s_l=n_lp_{l-1}+p_{l-2}-3-(n_l-1)
p_{l-1}=p_{l-1}-3+p_{l-2}$. Now ($5_{j-1}$) with $d=1$ implies ($4_j$) provided 
that $s_l=p_{l-1}(n_l-1)$, and in this case equality in ($4_j$) requires that 
$s_1=n_1-2$ and $s_l=p_{l-1}n_l$ for $l=2,\dots,j-1$.

Further, if $s_l>p_{l-1}(n_l-1)$, then $\lc s_l/p_{l-1}\rc\ge n_l$, while the 
inequality $s_1+\lc s_2/p_1\rc+\ldots+\lc s_{l-1}/p_{l-2}\rc\ge n_1+\ldots+
n_{l-1}-2$ follows from ($3_{l-1}$), ($5_{l-1}$), and statement d) for $j=l-1$. 
Summing the last two inequalities, we get $s_1+\lc s_2/p_1\rc+\ldots+\lc s_l/
p_{l-1}\rc\ge n_1+\ldots+n_l-2$, which is stronger than~($4_l$).

Finally, suppose that $s_l<p_{l-1}(n_l-1)$. Define a positive integer $d$ by
$\lc s_l/p_{l-1}\rc=n_l-d$. If $d=1$, we still have  $\lc s_l/p-{l-1}\rc=n_l
-1$, but $s_1+\ldots+s_{l-1}=p_l-3-s_l>p_{l-1}-3+p_{l-2}$, and the inequality 
($5_{j-1}$) with $d=1$, combined with the statement d) for $j=l-1$, implies 
that $s_1+\lc s_2/p_1\rc+\ldots+\lc s_l/p_{l-1}\rc>n_1+\ldots+n_{l-1}-2+n_l-1$, 
which is stronger than~($4_l$). If $d>1$, then $s_1+\ldots+s_{l-1}=p_l-3-s_l\ge
dp_{l-1}-3+p_{l-2}>p_{l-1}-3+dp_{l-2}$, and then $s_1+\lc s_2/p_1\rc+\ldots+\lc 
s_{l-1}/p_{l-2}\rc>n_1+\ldots+n_{j-1}-3+d$. Together with $\lc s_l/p_{l-1}\rc=
n_l-d$, this yields $s_1+\lc s_2/p_1\rc+\ldots+\lc s_l/p_{l-1}\rc>n_1+\ldots+
n_{l-1}-2+n_l-1$, which is stronger than~($4_l$). This proves statements a) and 
b) for~$j=l$. The proof of the statements c) and d) for $j=l$ is similar. 
Lemma~9 is proven.						\qed\enddemo

By Lemma~7, a rotation of~$\Delta_1$ by $2\pi q/p$ repeated $p_{i-1}$ times 
(where $i=1,\dots,k$), being equivalent to a rotation by $\pm2\pi r_i/p$, 
destroys all diagonals of $\Delta_1$ that are longer than~$r_i$. Repeating the 
argument from the case $k=2$ above (the case of a two-term continued fraction), 
we conclude that a rotation of~$\Delta_1$ by $2\pi q/p$ destroys at least $\lc 
s_i/p_{i-1}\rc$ diagonals in any length group $S_i$. Therefore, it destroys at 
least $s_1+\lc s_2/p_1\rc+\ldots+\lc s_k/p_{k-1}\rc$ diagonals of~$\Delta_1$. 
Setting $j=k$ in Lemma~9\,a), we get the statement of Theorem~4, provided that 
$k>1$, which is equivalent to $q\ne1$. The case $q=1$ was considered above. The 
proof of Theorem~4 is completed. 				\qed

\smallskip

One can show (see~\cite2 and \S7 below) that the lower bound given by Theorem~4 
is exact. The proof of Theorem~4 shows that the numbers of diagonals of the 
triangulation in the length group $S_i$ in this case are $s_1=n_1-2$, $s_i=p_
{i-1}n_i$ for $i=2,\dots,k-1$, and $s_k=p_{k-1}(n_k-1)$ (provided that $k>1$), 
which are the numbers given by Lemma~9\,b. Furthermore, the rotation by $2\pi 
q/p$ destroys only $(1/p_{i-1})$th fraction of the $s_i$ diagonals belonging to 
$S_i$, and the contributions to $E(p,q)-3$ from the groups $S_1,\dots,S_k$ are 
$n_1-2,n_2,\dots,n_{k-1},n_k-1$ respectively. 

\head\S7. Optimal triangulations of the regular $p$-gon \endhead

We show in this section that the bounds of Theorems~3 and~4 are sharp.

\proclaim{Theorem~5} For any $x\in\lpq\sm\{S^1_w,S^1_z\}$, the cut locus $C(x)$ 
is a simple spine of $\lpq$ with $E(p,q)-3$ vertices, provided that $q>1$. 
\endproclaim

\demo{Proof} This follows from the description of the convex hull of the 
$p$-element preimage of $x$ under the universal covering $S^3\to\lpq$, 
see~\cite2, and from the results of~\S5 and~\S6 above. 		\qed\enddemo

\proclaim{Theorem~6} Let $p$ and $q$ be coprime positive integers, $p\ge3$. 
\roster
\item"a)" There exists a triangulation $\Delta_1$ of a regular $p$-gon that can 
be converted into $\Delta_2$ by just $E(p,q)-3$ flips, where $\Delta_2$ is 
$\Delta_1$ rotated by $2\pi q/p$. 
\item"b)" Consider the cut locus $C(x)\subset\lpq$ of an arbitrary point $x\in
\lpq$ with respect to the standard metric on~$\lpq$. There exists a small 
perturbation $P$ of $C(x)$ such that $P$ is a simple spine of~$\lpq$ with 
exactly $E(p,q)-3$ vertices.
\endroster
\endproclaim

\demo{Proof} As in the proof of Theorem~4, we assume that $q<p/2$. The case $p=
3$ is trivial. Below we assume that $p\ge4$.

First consider the case $q=1$. By $A_1,A_2,\dots,A_p$ denote the consecutive 
vertices of the $p$-gon. Let $\Delta_1$ consist of the diagonals $A_1A_i$, 
where $i=3,4,\dots,n-1$; then $\Delta_2$ consists of the diagonals $A_2A_i$, 
where $i=4,\dots,n$. Take $\Delta_1$ and make a flip along $A_1A_3$ (for any 
diagonal, a unique flip involving it is possible), then along $A_1A_4$, and so 
forth. After the first flip we get $A_2A_4$, after the second one, $A_2A_5$, 
and so forth. In $p-3$ flips we get~$\Delta_2$. This proves statement a) for 
$q=1$. 

If $q=1$, then the combinatorial type of $C(x)\subset\lpq$ does not depend 
on~$x$; moreover, there exists an isometry of $L_{p,1}$ taking $x_1$ to $x_2$ 
and $C(x_1)$ to $C(x_2)$ for any $x_1,x_2\in L_{p,1}$. The 1-dimen\-sional 
stratum of $C(x)$ is a closed line $S^1$ of multiplicity~$p$, and a section of 
$C(x)$ transversal to $S^1$, denoted by $Y_p$, gets rotated by $2\pi/q$ by 
monodromy along~$S^1$, see~\S5. A sequence of $p-3$ flips converting $\Delta_1$ 
to $\Delta_2$ corresponds to a simple spine $P$ of $L_{p,1}$ with $p-3$ 
vertices, see the proof of Theorem~3. Thus, statement b) of the theorem in the 
case $q=1$ follows from the statement~a). 

Now suppose that $q>1$. If $x\in\lpq$ does not lie on the core circles $S^1_w,
S^1_z$ (see~\S5), then the statement~b) follows from Theorem~5; a 
perturbation of $C(x)$ is not necessary at all in this case. If $x$ lies on a 
core circle, it is enough to shift $x$ slightly out of $S^1_w$ or $S^1_z$. The 
statement~a) is deduced from the statement~b) by repeating a part of the proof 
of Theorem~3. Thus, Theorem~6 follows from Theorem~5.	\qed\enddemo

\midinsert
\epsfxsize=200pt
\centerline{\epsffile{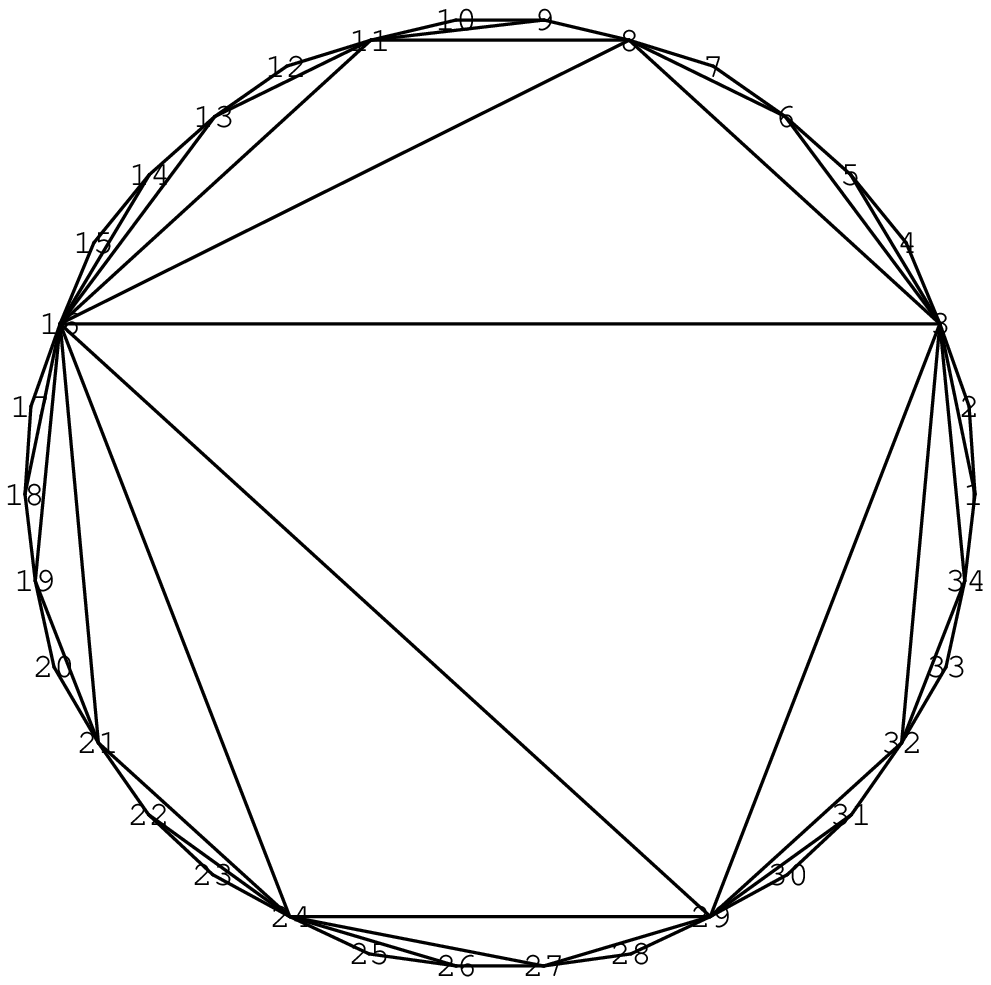}}
\botcaption{Figure 6} Optimal triangulation for $(p,q)=(34,13)$\endcaption
\endinsert

The following Mathematica${}^\circledR$ script\footnote{Special thanks to 
Roderik Lindenbergh and Martijn van Manen}\ produces the optimal (in the sense 
of Theorem~6) triangulation $\Delta_1$ for the case $p=34$, $q=13$ (then $k=6$, 
$n_1=2$, $n_2=n_3=n_4=n_5=1$, $n_6=2$, and $E(p,q)-3=5$). The best choice of 
the eccentricity {\tt e} depends on the pair~$(p,q)$.
\smallskip

\noindent{\tt <<DiscreteMath`ComputationalGeometry`}
\smallskip

\noindent{\tt e=.007; p=34; q=13;}\newline
{\tt points=Table[(1+e Sin[(2k Pi)(q/p)])}\newline
$\phantom{-----}${\tt \{Cos[(2k Pi)/p],Sin[(2k Pi)/p]\},\{k,0,p-1\}];}\newline
{\tt d=DiagramPlot[points, TrimPoints->3];}\newline
{\tt triangulation=PlanarGraphPlot[points];}\newline
{\tt Show[\{d,Graphics[RGBColor[1,0,0]],triangulation\},}\newline
$\phantom{-----}${\tt PlotRange->\{\{-1.2,1.2\},\{-1.2,1.2\}\}];}

\smallskip\noindent
The output is shown on Fig.~6.

\comment
$$\frac pq=n_1+\cfrac1\\
	n_2+\cfrac1\\
	\ddots\lower6pt\hbox{$+\dfrac1{n_k}$}\endcfrac .$$
\endcomment

\head Appendix: properties of the Farey tesselation \endhead

\demo{Proof of Lemma~\rom1} The modular group has two generators, $S\:z\mapsto
-\frac1z$ (central symmetry of~$H^2$ at~$i$) and $T\:z\mapsto z+1$, see~\cite8.
By construction of the Farey tesselation, both $S$ and~$T$ preserve it. This 
implies statement~c).

Statement b) is proved by induction over $E(p,q)$. It holds if $E(p,q)\le1$ (in 
this case $p/q=0$ or $p/q=\infty$ or $p/q=\pm1$). Let $E(p,q)>1$. Since the 
symmetries in the lines $(0,\infty)$ and $(1,-1)$ preserve the Farey 
tesselation, we can suppose that $p>q>0$. By the induction hypothesis, a 
geodesic segment $(ti,(p-q)/q)$ intersects (the interiors of) $E(p-q,q)$ Farey 
triangles. The isometry $z\mapsto z+1$ preserves the tesselation and takes $ti$ 
to $1+ti$ and $(p-q)/q$ to $p/q$, so that a geodesic line $(ti,(p-q)/q)$ gets 
mapped into a geodesic line $(1+ti,p/q)$ intersecting the same number $E(p-q,q)$
of Farey triangles. Since $p>q>0$, a geodesic line $(ti,p/q)$ intersects the 
Farey edge $(1,\infty)$ at $1+t'i$ (where ${t'}^2=p/q+1+t^2(1-q/p)$). The 
segment $(ti,1+t'i)$ of this line intersects the triangle $(0,1,\infty)$ and 
the remaining part of it intersects $E(p,q)-1$ other Farey triangles, so the 
line $(ti,p/q)$ cuts $E(p,q)$ Farey triangles. This applies also for the 
geodesic line $(0,p/q)$. Statement~b) follows. 

To prove statement a), note that the modular group action actually defines the 
Farey tesselation, not only preserves it. Indeed, maps $S$, $T$, and $ST^{-1}S$ 
take the triangle $(0,1,\infty)$ to its images under reflections in its sides. 
Now, if a triangle~$\Delta$ is obtained from $(0,1,\infty)$ by a modular 
transformation~$F$, one gets its three neighbors from $(0,1,\infty)$ by a 
modular transformation~$F^{-1}GF$ with $G=S$ or $T$ or~$ST^{-1}S$. Further, 
modular maps do not change $mq-np$, because $m'q'-n'p'=\det\pmatrix m'&p'\\ n'&
q'\endpmatrix=\det\(\pmatrix a&b\\ c&d\endpmatrix\pmatrix m&p\\ n&q\endpmatrix\)
=\det\pmatrix a&b\\ c&d\endpmatrix\det\pmatrix m&p\\ n&q\endpmatrix=mq-np$: 
since the matrix $\pmatrix a&b\\ c&d\endpmatrix$ of a modular map has 
determinant~1, it keeps $m'$ coprime with~$n'$ and $p'$ coprime with~$q'$.

Thus any Farey edge can be mapped to an edge of the triangle $(0,1,\infty)$ by 
a modular map, which preserves $mq-np$, but this amounts to $\pm1$ for any edge
of~$(0,1,\infty)$. This proves the ``only if'' part of statement~a). A mapping 
$z\mapsto\frac{mz+p}{nz+q}$, which is modular if $mq-np=1$, takes an edge $(0,
\infty)$ to $(m/n,p/q)$, which is a Farey edge as well by statement~c). If $mq-
np=-1$, consider the map $z\mapsto\frac{pz+m}{qz+n}$. Statement~a) follows.

In the assumptions of statement d) we have $mq-np=\pm1$ by virtue of~a). Then 
$m(n+q)-n(m+p)=\pm1$, and statement a) implies that $\(\frac{m+p}{n+q},\frac mn
\)$ is an edge of the Farey tesselation. Similarly, $\(\frac{m+p}{n+q},\frac pq
\)$, $\(\frac{m-p}{n-q},\frac mn\)$, and $\(\frac{m-p}{n-q},\frac pq\)$ are 
edges of the Farey tesselation, which thus contains triangles $\(\frac mn,\frac 
pq, \frac{m+p}{n+q}\)$ and $\(\frac mn,\frac pq, \frac{m-p}{n-q}\)$; this 
proves statement~d).

Statement e) is an immediate consequence of~d), by construction of the Farey 
series. By the way, this explains the term ``Farey tesselation''. 
\qed\enddemo

\enddocument
\end